\documentclass{amsart}
\usepackage{amsthm}
\usepackage{amsmath}
\usepackage{amsfonts}

\title{Counting Solutions to Binomial Complete Intersections}

\author{Eduardo Cattani and Alicia Dickenstein}
\address{Eduardo Cattani: Department of Mathematics
and Statistics. University
of Massachusetts. Amherst, MA 01003, USA}
\email{cattani@math.umass.edu}
\thanks{E. Cattani was partially supported by NSF Grant
DMS--0099707.  Part of this work was done while he was visiting
the University of Buenos Aires supported by a Fulbright Fellowship
for Lecturing and Research; he is grateful for their hospitality and
sponsorship.}

\address{Alicia Dickenstein: Departamento~de
Matematica, FCEyN.
Universidad de Buenos Aires. (1428) Buenos Aires,
Argentina}
\email{alidick@dm.uba.ar}
\thanks{A. Dickenstein is partially supported  by
UBACYT and CONICET, Argentina.}

\newcommand{\baseRing}[1]{\ensuremath{\mathbb{#1}}}
\newcommand{\Z}{\baseRing{Z}}
\newcommand{\R}{\baseRing{R}}

\newcommand{\N}{\baseRing{N}}
\newcommand{\Q}{\baseRing{Q}}

\theoremstyle{plain}
\newtheorem{theorem}{Theorem}[section]
\newtheorem{lemma}[theorem]{Lemma}
\newtheorem{corollary}[theorem]{Corollary}

\newtheorem{proposition}[theorem]{Proposition}

\theoremstyle{definition}

\newtheorem{definition}[theorem]{Definition}
\newtheorem{remark}[theorem]{Remark}

\newtheorem{example}[theorem]{Example}

\numberwithin{equation}{section}

\DeclareMathOperator{\supp}{supp}

\newcommand{\Script}[1]{\ensuremath{{\mathcal{#1}}}}

\newcommand{\UU}{\Script{U}}

\newcommand{\FF}{\Script{F}}

\newcommand{\JJ}{\Script{J}}
\newcommand{\LL}{\Script{L}}
\newcommand{\RR}{\Script{R}}
\newcommand{\VV}{\mathbb{V}}

\newcommand{\qq}{\mathfrak{q}}
\newcommand{\pp}{\mathfrak{p}}
\newcommand{\mm}{\mathfrak{m}}
\newcommand{\ee}{\mathfrak{e}}

\newcommand{\p}{p_1,\dots,p_n}
\newcommand{\tildep}{\tilde p_1,\dots,\tilde p_{n-1}}
\newcommand{\pn}{p_1,\dots,p_n}
\newcommand{\pc}{p_1(c;x), \dots, p_n(c,x)}

\newcommand{\kn}{k[x_1,\dots,x_n]}

\newcommand{\bk}{\bar k}

\newcommand{\M}{\mathcal{M}}

\renewcommand{\supp}{\rm supp}

\begin{document}

\begin{abstract}
We study the problem of counting the total
number of  affine solutions 
of a system of $n$ binomials in $n$ variables over
an algebraically closed field of characteristic zero.
We show that we may decide in polynomial time if that
number is finite.  We give a 
combinatorial formula for computing the total number of affine
solutions (with or without multiplicity) from which we deduce that
this counting problem is $\#P$-complete.  We discuss special
cases in which this formula may be computed in polynomial time; 
in particular, this is true for generic exponent vectors.
\end{abstract}
%\footnotetext[1]{AMS Subject Classification:
%Primary 13P10, Secondary 68Q17,
%14M25, 13P99. Keywords: binomial system,  counting problem, 
%complete intersection, $\#P$-complete.}

\maketitle

\section {Introduction}

A {\em binomial ideal} in the ring $\kn$ of polynomials with coefficients in
a field $k$, is an
ideal generated by binomials: $a x^\alpha - b x^\beta$, where
$\alpha, \beta \in \N^n$ and $a,b \in k^*$.  Binomial ideals
are quite ubiquitous in very different contexts particularly
those involving toric geometry and its applications
\cite{ eisenbud, sturmfels},
in the study of semigroup algebras, and in the modern versions
of hypergeometric systems of differential equations \cite{sst, dms}.
While binomial ideals are quite amenable 
to Gr\"obner and standard
bases techniques \cite{km1, km2}, they also provide some of the 
``worst-case"
examples in computational algebra, such as the Mayr-Meyer ideals
\cite{mayr}.

In this paper we consider ideals generated by $n$ binomials
in $R := \kn$, with ${\rm char}(k)=0$.  
Let $\bar k$ denote the algebraic closure of $k$.  We are interested in 
determining when  the number of solutions in $\bar k^n$ 
is finite and non zero (i.e., when the given binomials  define
a complete intersection in $R$) and, in this case, to count the number
of solutions, with or without multiplicity.
We will obtain 
properties of  these ideals
directly in terms of the given data: the exponents $\alpha, \beta$,
and the coefficients $a, b$.

Our starting point is then  a system of $n$ binomials in $R$, with
non-zero coefficients. Thus,
 we may assume that they are of the form
 \begin{equation}\label{thesystem}
 p_j(c;x) \ :=\  x^{\alpha_j} - c_j x^{\beta_j}\,;\quad  j=1,\dots,n,
\end{equation}
where
$ \alpha_j, \beta_j \in \N^n, \  \alpha_j \not= \beta_j$.
Let  $\JJ$ be the ideal generated by $p_1,\dots,p_n$ in the
polynomial ring $k(c)[x]$.
Given a choice of coefficients $c\in (k^*)^n$, let $\JJ_c$ be the ideal in $R$ generated
by $p_1(c;x), \dots, p_n(c,x)$ and $\VV_c \subset \bk^n$ the variety 
defined by $\JJ_c$.

Proposition~\ref{prop:soltorus}, which is a restatement of 
results in \cite{eisenbud}, gives a complete picture of the
number of solutions of the system (\ref{thesystem}) in the
algebraic torus $(\bar k^*)^n$.  Let $B$ be the matrix
\begin{equation}\label{matrixb}
B\ :=\ \left(
\begin{array}{c}
\alpha_1 - \beta_1\\
\alpha_2 - \beta_2\\
\vdots\\
\alpha_n - \beta_n
\end{array}\right)\,,
\end{equation}
whose
$j$-th row is the vector $\alpha_j - \beta_j$. 
Then, for generic coefficients $c \in (k^*)^n$, 
$\VV_c \cap (\bk^*)^n$ consists of $|\det B|$-many points
all of which have multiplicity one
(this may be seen directly or as a simple instance of 
Bernstein's theorem).  
In fact,
 if 
$\det B \not =0$,  this is true for all $c\in (k^*)^n$.  On the other hand,
if $\det B  =0$, then, for  coefficients $c\in (k^*)^n$ not satisfying
 the algebraic conditions (\ref{citorus})  it holds that $\VV_c \cap (\bk^*)^n = \emptyset$ ,
while if the coefficients
satisfy  (\ref{citorus}), the variety
$\VV_c \cap (\bk^*)^n$ has codimension equal to the rank of $B$.
We set $\delta:=|\det B|$.

Deciding whether the system (\ref{thesystem}) has a non-empty, finite
set of solutions in $\bk^n$ is more involved.  We must,
first of all, consider the possibility that some 
exponent vector $\alpha_j$ or $\beta_j$ may vanish.  This is equivalent
to the statement that some variables $x_j$ are invertible modulo the
ideal $\JJ$.  The reduction to the case when this does not happen
is accomplished in Proposition~\ref{invertible}.  We may then assume
that $0\in \VV_c$ for all choice of coefficients.  Now, in the generic
case $\det B \not=0$, Theorem~\ref{fscriterion} gives a condition on 
the exponents of the system that guarantees that the system (\ref{thesystem})
is a complete intersection for all $c\in (k^*)^n$.  If, on the other hand,
$\det B =0$, Theorem~\ref{fscriterion} only implies that 
(\ref{thesystem})
is a complete intersection for a generic set of coefficients $c\in (k^*)^n$.
Indeed, in this case, algebraic conditions such as (\ref{citorus})  enter
into play.  This leads to the notion of {\em generic complete intersection}, 
that we will abbreviate by {\em  gci}.
We will say that $\pn$ is a gci if  $\JJ_c$ 
is a complete intersection  in
$R $, i.e., $\VV_c$ is a finite non empty set, for generic coefficients $c\in (k^*)^n$.  

Even though Theorem~\ref{fscriterion} gives a combinatorial criterion
for deciding if $\pn$ is a gci, its verification requires $2^n$ steps.  
One of the main results
of this paper is Theorem~\ref{ciprocedure} where we describe
a polynomial-time algorithm  to decide  whether 
$\pn$ is a gci directly from the exponents $\alpha_j,\beta_j$.

Given a generic complete intersection $\pn$, let 
\begin{equation}\label{dimension}
d\ :=\ \dim_k \kn/\JJ_c\ ;\quad D \ :=\ \dim_k \kn/\sqrt{\JJ_c}
\end{equation}
be the total number of points 
in
the variety $\VV_c$, counted with and without multiplicity.  Given an index set
 $L \subset \{1,\dots,n\}$,  we denote by $\mu_L$, the number of points in $\VV(\JJ) \cap \bk^n_L$,
$\ \bk^n_L \ :=\ \{x\in \bk^n : x_\ell=0 \ \hbox{if and only if}\ \ell\in 
L\}\,
$, counted with multiplicity.
  We set $[n]: = \{1, \dots, n\}$ and $\mu := \mu_{[n]}$, the multiplicity
 at the origin.

In Section~\ref{sec:irreducible} we compute $d$, $D$, and $\mu_L$ for a
gci.  A key 
ingredient
is what we call {\em parametric reduction}, which allows us to reduce 
the study of generic
complete intersection binomial ideals to a particular class of ideals with a
normalized presentation.  We show in
Theorem~\ref{parred} that we can keep track of
the various multiplicities through the process
of parametric reduction. We then compute $d$ and $D$ for so-called
{\em irreducible systems}.  We show that  an irreducible
system that is in normal form may behave in one of three possible
ways: its binomials are a standard
basis for either a global or a local term order,
or they are weighted homogeneous.  This allows us to
read off  the dimension and multiplicities from the exponents (cf. 
Theorem~\ref{vinberg}).
 Interestingly,
the linear algebra problem that underlies these results 
appeared in the work of Vinberg about Cartan matrices
\cite[Theorem 4.3]{kac}.  
For generic exponents, a binomial system in normal form is irreducible
and has $\det B\not=0$.
Hence, Theorem~\ref{vinberg} gives a 
polynomial time algorithm for computing the number of solutions
of a  complete intersection binomial system
with generic exponents and arbitrary non-zero coefficients.

We next consider the case of a general gci.  
Using a well-known quadratic-time algorithm, due to Tarjan \cite{tarjan},
we find a block decomposition of the system into irreducible ones.
{} From this
decomposition we construct an acyclic directed graph naturally
attached to the system.  In Theorem~\ref{th:dformula} we give an
explicit combinatorial formula to compute the
dimensions and multiplicities of the system from this graph.

Section~\ref{sec:complexity} is devoted to  counting complexity issues.
We reverse the correspondence from binomial systems
to acyclic digraphs and assign to each such graph a simple binomial
system. The number of solutions of this system corresponds to  
invariants of the
graph whose computation is known to be $\#P$-complete.
Indeed, we show that  particular instances correspond to
counting independent sets in  bipartite graphs, or more generally,
antichains in a poset; both of these problems are known to be 
$\#P$-complete \cite{valiant, provanball}.  
Hence, even though the problem of deciding whether a system is a 
gci as well as the problem of counting the number
zeros in the torus of the binomial system
 defined by (\ref{thesystem}), are solvable in polynomial time, 
we prove in Theorem~\ref{th:sharpP} that counting the total number
of affine solutions, with or without multiplicity, is a
$\#P$-complete problem.  Thus, binomial systems furnish a very
simple example of the type of problems, ``easy" to decide but 
``hard" to count that motivated 
Valiant's introduction of the notion of counting complexity \cite{valiant}.  
Finally, in Proposition~\ref{polycase} we
identify another class of systems whose solutions may be computed 
in polynomial time.

The last section of the paper is devoted to a brief discussions of some of the 
applications of this work which motivated our study.
We show, first of all, how 
Theorem~\ref{th:dformula} may be applied      to compute the
multiplicity and geometric degree
\cite{bm} of the primary components of a lattice basis
ideal 
$J\subset k[x_1,\dots,x_m]$.  
This, in turn, may be used to describe the holonomic
rank of Horn systems of hypergeometric partial differential equations and to
study sparse discriminants,
generalizing the  codimension-two case. \cite{ ds, dms}.
Finally we recall the results of \cite[Chapter 10]{stubook}
relating the study of systems of partial differential equations
with constant coefficients with that of the corresponding algebraic
system. 

\section{Complete Intersections and normal forms}
\label{sec:compint}

We begin by considering the question of when  binomials
$\pc$ as in (\ref{thesystem})
  define a complete intersection when viewed as
elements of the Laurent polynomial ring
$S := k[x_1^{\pm 1},\dots, x_n^{\pm 1}]$.
Let $B$ be the $n\times n$ exponent matrix 
defined in (\ref{matrixb}).
 We note that
even though the
rows of $B$ are only defined up to sign, this will not affect
our arguments.
It follows from
\cite[Theorem~2.1]{eisenbud} that if $\det B\not=0$
then,
for any choice of coefficients in $(k^*)^n$, 
$p_1(c;x),\dots,p_n(c;x)$ define a regular sequence in $S$.
Moreover,
the system of equations
\begin{equation}\label{eqn}
p_j(c;x)\ =\ 0\ ;\quad j=1,\dots,n
\end{equation}
has $|\det B|$-many solutions in the algebraic torus
$(\bk^*)^n$  and all of them are
simple.

On the other hand, if $\det B =0$ then $p_1(c;x),\dots,p_n(c;x)$ does 
not define
a complete intersection in $S$
for any choice of coefficients.
Indeed, if the system
(\ref{eqn}) has a solution  $ x \in (\bk^*)^n$,
it will necessarily have infinitely many.
Let $\RR$ be the lattice of relations
$$\RR \ :=\ \{m\in \Z^n : \sum_{j=1}^n m_j(\alpha_j - \beta_j) = 0\}.$$
For any $m\in \RR$ we have a $\bk^*$-action on the set of solutions
of (\ref{eqn}) defined by $(t;x) \mapsto (t^{m_1} x_1,\dots,t^{m_n} 
x_n)$,
and therefore the set of solutions could never be finite.
Note also that if $\det B = 0$ then, for generic coefficients $c_j$,
(\ref{eqn}) has no solutions. In fact,
if $x\in (\bk^*)^n$ is a solution of (\ref{eqn})
we have
$$x^{\alpha_j - \beta_j} = c_j, \, {\rm for \, all\/} \, \, j \, = \, 1, \dots, n,$$
and therefore  
$$\prod_{j=1}^n\,c_j^{m_j} \ =\ 1\,, \, {\rm for \, all\/}\,  \, m\in \RR.
$$
Thus, if $\nu^1,\dots,\nu^r$ is a basis of $\RR$, a necessary condition
for $\pc$ to have a solution in $(\bk^*)^n$ is that,
\begin{equation}
\label{citorus}
\prod_{j=1}^n\,c_j^{\nu^\ell_j} \ =\ 1\ \hbox{ for all  \, 
$\ell = 1,\dots,r$}.
\end{equation}
This condition is also sufficient.
Suppose that (\ref{citorus}) holds and let
$\LL$ be the sublattice of $\Z^n$ spanned by $\alpha_j - \beta_j$,
$j=1,\dots,n$. Denote by $\rho \colon \LL \to \bk^*$ the
group homomorphism
(i.e., the partial character) defined by
$$\rho(\alpha_j - \beta_j) = c_j\,.$$
The equalities in (\ref{citorus}) imply that $\rho$ is well-defined
and, since up to a monomial (which is invertible
in the Laurent polynomial ring),
$$ p_j(x) \, = \, x^{\alpha_j - \beta_j} - \rho({\alpha_j - \beta_j})$$
it follows from \cite[Theorem~2.6]{eisenbud} that $\pc$
define an ideal in $S$ of codimension equal to the rank of $\LL$.
Hence we obtain:

\begin{proposition}\label{prop:soltorus}
Let $\pc$ be as in (\ref{thesystem})
and $B$ as above.  For any choice of coefficients
$c \in (k^*)^n$, the ideal  they generated 
in $S$  is 
a complete intersection  if and only if
  $\det B  \not=0$.  If  $\det B =0$ and the identities (\ref{citorus}) are satisfied
  then the binomials (\ref{thesystem}) define an ideal in
  $S$ of codimension equal to the rank of $B$.
\end{proposition}

In the remaining part of this section, we
will discuss criteria for deciding when $\pn$ is a gci.
Since we
are not assuming that ${\rm supp}(\alpha_j) \cap {\rm supp}(\beta_j) = \emptyset$, 
where, for $v\in \R^n$:
$${\rm supp}(v) := \{i\in [n] : v_i \not =0\},$$
the matrix $B$, by itself, does not
allow us to recover the exponents of the  binomials (\ref{thesystem}).
It is useful to introduce the following concept, already present
in the work of Scheja, Scheja, and Storch \cite{sss}:

\begin{definition} \label{def:ZofK}
Let $p_j =\, x^{\alpha_j} - c_j \, x^{\beta_j}$, $j=1,\dots,n$, be
a system of binomials in $k[x_1,\dots,x_n]$.  For each index set
$K \subset [n]$, let
\begin{equation}\label{zeroset}
Z(K) \ :=\ \{j\in [n] : \supp(\alpha_j) \cap K \not=\emptyset
\ \hbox{and}\ \supp(\beta_j) \cap K \not=\emptyset\}\,.
\end{equation}
\end{definition}

We start by showing that we can restrict ourselves to the case where
$0\in \VV_c$.  Since this property is equivalent 
to the statement
that all exponent vectors are non zero, it is independent of the
choice of coefficients.  We want to identify all indices $i$
for which $x_i$ is invertible modulo the ideal $\JJ$, i.e., the $x_i$ coordinate of any solution to the system of
binomials is necessarily non zero.
Set $I_0 = \emptyset$ and, for $\ell \geq 1$, let 
$$I_\ell := \bigcup_{} \{{\rm supp}(\alpha_j):{\rm supp}
(\beta_j)\subset I_{\ell-1}\} \cup \bigcup_{} \{{\rm supp}(\beta_j):{\rm supp}(\alpha_j)\subset I_{\ell-1}\}$$
and $I = \bigcup_\ell I_\ell$.  Induction on $\ell$ shows
easily that if $i\in I$, the variable $x_i$ is invertible
modulo the ideal $\JJ$ and, conversely, that these are all the
variables invertible modulo $\JJ$.  Thus, after
reordering of variables and polynomials, we may assume
that the variables $x_{r+1},\dots,x_n$ are invertible and
that the binomials $p_{s+1},\dots,p_n$ involve only
the variables $x_{r+1},\dots,x_n$, while for $j\leq s$
both monomials $x^{\alpha_j}$ and $x^{\beta_j}$ are divisible by
at least one of the variables $x_i$, $i\leq r$, i.e., that $Z([r])
=  [s]$.
Following \cite{fs} we define:

\begin{definition}\label{derived}
Let $x':=(x_1,\dots,x_r)$, $c' := (c_1,\dots,c_r)$.
For $j\leq s$, set
\begin{equation}\label{phat}
\hat p_j(c';x') = p_j(c';(x_1,\dots,x_r, 1,\dots,1))\,.
\end{equation}
Then, the binomial system $\{\hat p_1,\dots,\hat p_s\}\subset
k(c')[x']$ is called the {\em derived system} of
$\pn$.  We denote by $\hat B$ the associated $s\times r$~matrix
as in (\ref{matrixb}).
\end{definition}

Note that $0 \in \VV(\hat p_1,\dots,\hat p_s)$ and 
that the matrix $B$ is of the form
$$B \ =\ 
\left(
\begin{array}{cc}
\hat B & * \\
0 & B_2
\end{array}
\right).
$$ 

\begin{lemma}\label{preinv}
Assume $\pn$ as in (\ref{thesystem}) is a gci and let $r, s$ be as above. Then,
$r=s$  and  $\det(B_2) \not=0$.
\end{lemma}

\begin{proof}
Since the variables $x_{r+1},\dots,x_n$ are all invertible 
modulo $\JJ$, the system of equations $p_{s+1} = \cdots = p_n =0$,
is
equivalent to the system $x^{\alpha_{j} -\beta_{j} } = c_j, \ 
\hbox{for all\ } j = s+1, \dots,n$.
Hence, arguing as in the discussion leading to
Proposition~\ref{prop:soltorus}, we see that each integer relation
among the vectors $\alpha_j - \beta_j, j = s+1, \dots, n$ imposes a
polynomial condition on the coefficients as in (\ref{citorus}). 
If $s < r$, then  $n-r < n-s$  and so there exists a non trivial
relation. Therefore, $p_1, \dots, p_n$ 
has generically no solutions, a
contradiction.  
On the other hand,  if $s > r$,  or if $r=s$ and $\det(B_2)=0$,  
then, generically, the system 
$p_{s+1}(x_{r+1}, \dots, x_n) = \dots = p_n(x_{r+1}, \dots, n)=0$
has either no solutions or infinitely many in $(\bk^*)^{n-r}$. 
Since any solution of these equations may be extended to a solution
of (\ref{eqn}) by setting $x_1 =\cdots=x_r=0$, 
we get a contradiction again.
So $s=r$ and $\det(B_2)\not=0$, as claimed. 
\end{proof}

\begin{proposition}\label{invertible}
Let $\pn$,  $B$ be as above. Assume that $s=r$ and $\det(B_2) \not=0$.
Let $\hat p_1,\dots,\hat p_r$ be the derived system.
Then
$\pn$ is a gci   if and only if $\hat p_1,\dots,\hat p_r$
is a gci.
\end{proposition}

\begin{proof}
Assume  $\pn$ is a gci and let $\mathcal U$ be an open dense subset of $(\bk^*)^{n}$
such that the binomials with coefficients in $\mathcal U$ define a complete intersection
ideal in $\bk[x_1,\dots, x_n]$. It suffices to show that
the intersection of $\mathcal U$ with the fiber
$(\bk^*)^r\times\{(1,\dots,1)\}$ is also Zariski dense in the fiber.
Let $a''\in (\bk^*)^{n-r}$ be such that 
$\UU \cap \left((\bk^*)^r\times\{a''\}\right)$ is Zariski dense. 
Let
$\lambda'' \in (\bk^*)^{n-r}$
be a common zero  of $p_{r+1}(a'';x), \dots, p_n(a'';x)$.  Then,
since $s=r$,
the change of variables that sends $x_i$ to itself for $i=1, \dots, r$ and
$$x_j \mapsto x_j/\lambda''_j, \quad j = r+1, \dots, n,$$
transforms any of the  last $n-r$ polynomials $p_j, j= r+1, \dots, n,$ into a
non-zero multiple of $x^{\alpha_j} - x^{\beta_j}$  and, for $i\leq r$, 
the binomial $p_i$  into a non-zero multiple of
$$ x^{\alpha_i} - (\lambda'')^{\alpha''_i - \beta''_i}\, c_i \, x^{\beta_i},$$
where $\alpha''_i, \beta''_i \in \N^{n-r}$ denote the vectors consisting
of the last $n-r$ coordinates of $\alpha_i, \beta_i$.
Since this scalar transformation in the coefficient space $(\bk^*)^r$ preserves Zariski dense
subsets our assertion follows.
\par
Conversely,  assume that $\hat p_1,\dots,\hat p_r$ is a gci
and that $\det(B_2)\not=0$. Let $\varphi$ be a
non zero polynomial such that $\varphi(c') \not =0$ for a given $r$-tuple of coefficients
$c' = (c_1,\dots, c_r)$ implies that the corresponding polynomials
$\hat p_1(c';x'), \dots, \hat p_r(c';x')$ 
define a complete intersection.  Denote as before $c'' =  (c_{r+1}, \dots, c_n)$
and consider the rational function
$$\psi(c', c'') =  \prod_{\lambda'' \in \VV_{c''}} 
\varphi ((\lambda'')^{\alpha''_1-\beta''_1}\,c_1,\dots,(\lambda'')^{\alpha''_r-\beta''_r}\,c_r).$$
If $\psi(c', c'')$ is defined and non zero, then for any choice of the $|\det (B_2)|$-many roots $\lambda''$ of
the last $n-r$ polynomials, the specialized system
$$p_1(c';(x',\lambda'')) = \cdots = p_r(c';(x',\lambda'')) =0$$
has finitely many solutions and, consequently,
 $\pn$ is a gci.
\end{proof}

The following result is a reformulation of Theorem~2.3 in \cite{fs}.

\begin{theorem}\label{fscriterion}
Let $\pn$ be as in (\ref{thesystem}) and
suppose that 
 $0 \in \VV(\JJ)$.
Then, $\pn$ ia a gci if and
only if $|Z(K)| \leq |K|$ for all $K \subset [n]$.
\end{theorem}

\begin{proof}
Suppose there exists $K\subset [n]$ such
that $|Z(K)| > |K|$.  Assume that $K$ is maximal with this
property.  After reordering, if necessary, we may assume 
that $K = \{r+1,\dots,n\}$ and $Z(K) = \{s+1,\dots,n\}$ where
$s<r$.  Since $0 \in \VV(\JJ)$, the maximality assumption implies
that the first $s$ binomials depend only on $x'=(x_1,\dots,x_r)$.
Otherwise, we may assume that there exists $k_1>r$, 
$k_1 \in {\rm supp}(\alpha_s)$.  Since $0 \in \VV(\JJ)$, 
$ {\rm supp}(\beta_s) \not= \emptyset$. If there exists
 $k_2>r$, 
$k_2 \in {\rm supp}(\beta_s)$, then $s \in Z(K)$ which is
a contradiction.  Therefore, 
${\rm supp}(\beta_s) \subset [r]$ and for any $\ell \in
{\rm supp}(\beta_s)$, $K' := K \cup \{\ell\}$ satisfies
$Z(K) \cup \{s\} \subset Z(K')$.  Hence $|Z(K')|>|K'|$ and
this contradicts the maximality of $K$.

Thus, for a given choice of coefficients, the system
\begin{equation}\label{e1}
p_1(c;x') = \cdots= p_s(c;x') = 0
\end{equation}
is either inconsistent or its solution space has dimension at
least $r-s>0$.  Since, any solution of (\ref{e1}) can be
extended to a solution of the full system by setting the $K$-coordinates
equal to zero, it follows that $\pn$ is not a gci.

Conversely, suppose $|Z(K)| \leq |K|$ for all $K \subset [n]$.
In order to show that $\pn$ is a gci it suffices to prove that
given any subset  $L\subset [n]$, for generic coefficients $\pc$
has at most finitely many solutions with zeros in $\bk^n_L$, where
\begin{equation}\label{zeroesL}
\bk^n_L \ =\ \{x\in \bk^n : x_\ell=0 \ \hbox{if and only if}\ \ell\in 
L\}\,.
\end{equation}
Assume that for some choice of coefficients, there exists a solution
in $\bk^n_L$.  Then, for any $i\not\in Z(L)$,
$p_i$ depends only on the variables in $J$, the complement of $L$ in $[n]$
and hence, since $0 \in \VV(\JJ)$, $Z(L)^c \subset Z(J)$.  Since, by
assumption $|Z(L)|\leq |L|$ and $|Z(J)|\leq |J|$, we deduce that
$$ |L| \leq |Z(J)^c|\leq |Z(L)| \leq |L|,$$
and therefore $|Z(L)|=|L|$.  Reordering we may assume that
$J=Z(J) = [r]$ and let $B_1(L)$ denote the $r\times r$ exponent matrix
as in (\ref{matrixb}).  If  $\det B_1(L) =0$, then for generic coefficients 
the first $r$ binomials have no solutions in $(\bk^*)^r$ and hence,
generically, $\pc$ have no solutions in $\bk^n_L$.
On the other hand,
if $\det B_1(L) \not=0$ then, for all choices
of coefficients in $(k^*)^r$, there exists finitely many solutions
of $p_1=\cdots=p_r=0$
in $(\bk^*)^r$ and hence finitely many solutions
of $\pc$ with zeros exactly in $L$.  
\end{proof}

\begin{remark}\label{setsL}
Note that in the proof of Theorem~\ref{fscriterion}
we have shown that if $\pn$ is a gci, $L\subset [n]$,
and $\bk^n_L$ is as in (\ref{zeroesL}),
then, for generic coefficients, there exists a solution in
$\bk^n_L$ if and only 
if $|Z(L)| = |L|$ and, after reordering so that $Z(L)=L=\{r+1,\dots,n\}$,
the binomials $p_1,\dots,p_r$ depend only on the first $r$ variables,
and the corresponding $r\times r$ exponent matrix
$B_1(L)$ is non-singular.  Moreover, 
for generic $c\in (k^*)^n$, there are $|\det B_1(L)|$-many 
points (counted without multiplicity) in $\VV_c\cap \bk^n_L$.  Then, the number
of points in $\VV_c$, counted without multiplicity, is given by
\begin{equation}\label{computingD}
D\ =\ \sum_{\mu_L\not=0} |\det B_1(L)|,
\end{equation}
where $\mu_L$ is the total number of points in $\VV_c\cap \bk^n_L$ counted
with multiplicity.
We will develop in Section \ref{sec:irreducible} the combinatorics needed to
describe all sets $L$ with $\mu_L \not=0$  and
 we shall show in Section \ref{sec:complexity} that counting the
number of such sets is a $\#P$-complete problem.
\end{remark}
Note that if 
$0\in \VV(\JJ)$, the condition
that $\pn$ is a gci 
depends only on the combinatorics
of the exponents $\alpha_j, \beta_j$.  
It follows from Proposition~\ref{prop:soltorus} and 
Theorem~\ref{fscriterion} than, when $\det(B) \not =0$,
if $\pn$ is a gci, then it is a complete intersection for
any choice of the coefficients (as long as $c_j\in k^*$).

The variant of the Fischer-Shapiro criterion
embodied in Theorem~\ref{fscriterion} allows us to determine
whether $\pn$ is a gci.  However, this
involves checking exponentially many conditions, one for
each subset $K \subset [n]$.  We will now  show how this can be done
in a number of steps that depends 
 polynomially (on $n$).  We begin with the following
simple corollary to Theorem~\ref{fscriterion}.

\begin{corollary}\label{cor1}
Suppose $\pn $ is a gci
and $0\in \VV(\JJ)$.  Let
$$\M \ = \ \{x^{\alpha_j}, x^{\beta_j}\ ;\ j=1,\dots,n\}$$ denote the
set of monomials appearing in $\pn$. Then for each $i\in [n]$ there 
exists
$r_i>0$ such that $x_i^{r_i} \in \M$.

\begin{proof}
If for some $i\in [n]$, $x_i^{r_i} \not\in \M$ for all $r_i > 0$, then $Z(\{1,\dots,\hat i,\dots,n\})= [n]$, contradicting
Theorem~\ref{fscriterion}.
\end{proof}

\end{corollary}

One can easily give examples showing that
the necessary condition in Corollary~\ref{cor1}
is not sufficient to guarantee that 
$\pn$ define a gci.
However, the following stronger notion provides a
sufficient condition.

\begin{definition} \label{def:nf}
We say that $\p$ are in {\em normal form} if and only if
for all $i\in [n]$
$$ p_i = x_i^{r_i} - c_i x^{\beta_i}\ ; \quad r_i>0, \, \beta_i\not=0.$$
Note that if the system is in normal form
then $0\in \VV(\JJ)$.
\end{definition}

\begin{proposition}\label{normalform}
Assume $\pn$ are in normal form. Then $\pn$ is a gci.
\end{proposition}
\begin{proof}
For any $K\subset[n]$, $Z(K)\subset K$ and the result follows from
Theorem~\ref{fscriterion}.
\end{proof}

We will next show how to reduce ourselves to systems
$p_1, \dots, p_n$ in normal form.

\subsection{Parametric Reduction}

Let $\pn$ be a binomial system and suppose that they satisfy
the necessary condition in Corollary~\ref{cor1}, but that it is not possible to relabel
variables and binomials, or invert the coefficient  of one or
more binomials, so as to put the system in normal form.
This means that one of the binomials must contain two monomials
of the form $x_i^{r_i}$ and $x_j^{r_j}$ with $i\not=j$.
Then, after relabeling 
  we may assume that $p_n$ is of the form
\begin{equation}\label{eq:parametric}
p_n  \ =\ x_n^{\ell} - c_n x_{n-1}^{m}  \,,\  \ell, m >0.
\end{equation}
Let $q := {\rm gcd}(m,\ell)$ and set
$m' := m/q$, $\ell' := \ell/q$.
We will consider the polynomial map that sends polynomials in
$n$ variables $x_1,\dots,x_n$ to polynomials in 
$n-1$ variables $u_1,\dots,u_{n-1}$:
\begin{equation}\label{def:parred}
x_i \mapsto u_i, \ i=1,\dots,n-2; \quad
x_{n-1} \mapsto u_{n-1}^{\ell'}; \quad
x_{n} \mapsto u_{n-1}^{m'}. 
\end{equation}
Let
$\tilde p_j$, $j=1,\dots,n-1$ be the image
of the  binomials $p_1,\dots,p_{n-1}$.
We will refer to $\tildep$ as a {\em parametric
reduction} of $\pn$ and denote by $\tilde\JJ$ the
ideal they generate in $k(c_1,\dots,c_{n-1})[u_1,\dots,u_{n-1}]$.

\begin{proposition}\label{parametricdet}
Suppose $\tildep$ is a parametric reduction of $\p$ and let $\tilde B$
and
$B$ be the associated matrices. Then $|\det B| = q\cdot |\det \tilde B| 
$.
Moreover, $0\in \VV(\JJ)$ if and only if  
$0\in \VV(\tilde\JJ)$ and, in this case, $\p$ is a gci if and 
only if $\tildep$ is
gci.
\end{proposition}

\begin{proof} The matrix $B$ is of the form
$$ B = \left(
\begin{array}{ccccc}
{} & {} & {} & {} & {}\\
\tilde b_1 &\dots & \tilde b_{n-2} & \tilde b_{n-1} & \tilde b_{n}\\
{} & {} & {} & {} & {}\\
0 & \dots & 0 & -m & \ell\\
\end{array}
\right)
$$
where $\tilde b_1,\dots,\tilde b_n$ are vectors in $ \Z^{n-1}$. On the 
other
hand, the matrix $\tilde B$ is given by
$$\tilde B = \left(
\begin{array}{cccc}
{} & {} & {} & {} \\
\tilde b_1 &\dots & \tilde b_{n-2} & \ell' \tilde b_{n-1} + m' \tilde 
b_{n}\\
{} & {} & {} & {}
\end{array}
\right)
$$
The first assertion now follows from a last-row expansion of $\det B$.

Suppose now
that $\p$ is not a gci.
By Theorem~\ref{fscriterion} there exists  $K\subset [n]$ such
that $|Z(K)|>|K|$. If $K \subset [n-1]$, then $Z(K)\subset [n-1]$
as well  and therefore by Theorem~\ref{fscriterion} $\tildep$ is not
a gci either. If $n\in K$, then taking $\tilde K =
K \backslash \{n\}$ we get that  $ Z(K) \backslash \{n\} \subset Z(\tilde K)$. Hence
 $|Z(\tilde K)| >
|\tilde K|$ and $\tildep$ is not a gci.

Conversely, if $\tildep$ is not a gci then there 
exists
$\tilde K \subset [n-1]$ such that $|Z(\tilde K)| > |\tilde K|$.
If $\tilde K \subset  [n-2]$ we take $K = \tilde K$ and then $Z(K) =
Z(\tilde K)$;
if, on the other hand, $n-1 \in \tilde K$, then we take
$K = \tilde K \cup \{n\}$ in which case $Z(K) = Z(\tilde K) \cup \{n\}$.
In either case $|Z(K)|>|K|$ and we are done.
\end{proof}

The results of this section may be summarized in a polynomial-time
algorithm to check whether a binomial system is a gci.

\begin{theorem}\label{ciprocedure}
We may decide in polynomial time 
whether  $\pn$ is a gci.  Moreover, if it is known
that $\det B\not =0$ we can check if $\pc$ is a
complete intersection in time $O(n^2)$.
\end{theorem}

\begin{proof}
It is easy to see from the procedure for constructing the
derived system that this step may be accomplished in 
at most $O(n^2)$ steps. If the number of non-invertible
variables does not equal the number of binomials in the
derived system then, by Lemma~\ref{preinv},
$\pn$ is not a gci.  Again by Lemma~\ref{preinv} we next check
whether $\det B_2\not=0$ (this is, of course, unnecessary
if it is known that $\det B\not=0$).  If so,
Proposition~\ref{invertible} allows us to restrict ourselves
to the derived system. We move down the list of binomials 
 searching for binomials of the form
$x_i^{r_i} - c x_j^{r_j}$.  Whenever such a binomial is
found we do parametric reduction and reduce by one the
number of binomials and of variables.  This step is then
repeated until there are no longer any binomials of that
form.  Clearly, this 
process stops after a quadratic number of steps.
Then $\pn$ is a gci if and only if Corollary~\ref{cor1}
holds.  This verification can certainly be carried
out in quadratically many steps.
\end{proof}

\begin{example}\label{exparred}
Consider the following binomials in $k[x_1,\dots,x_8]$:
\begin{eqnarray*}\label{ex:system}
p_1 = x_1^2 - x_2^3 ;\quad
p_2 = x_1 x_2 - x_1 x_3; \quad
p_3 = x_1^2 x_2 x_3 - x_3^7;\\
p_4 = x_4^2 - x_1^2 x_4^3;\quad
p_5 = x_5^2 - x_6^4;\quad
p_6 = x_5 x_6 - x_2 x_3 x_7^2 x_8;\\
p_7 = x_5 x_7 - x_7^2;\quad
p_8 = x_8^3 - x_1 x_6 x_7 x_8,
\end{eqnarray*}
where, since $\det B \not=0$,
we have set all coefficients $c_j =1$.
Although the system satisfies the necessary condition in
Corollary~\ref{cor1}, it is not in normal form.  We may
apply parametric reduction simultaneously to the binomials
$p_1$ and $p_5$ by considering the polynomial map from
$k[x_1,\dots,x_8]$ to $k[u_1,\dots,u_6]$
that sends:
\begin{eqnarray*}
x_1 \mapsto u_1^3;\quad x_2 \mapsto u_1^2;\quad x_3 \mapsto u_2;\quad 
x_4 \mapsto u_3;\\
x_5 \mapsto u_4^2;\quad x_6 \mapsto u_4;\quad x_7 \mapsto u_5;\quad 
x_8 \mapsto u_6.
\end{eqnarray*}
Here we have taken into account that the gcd of the exponents in
$p_5$ is $2$.  After changing signs when necessary, 
the new system $\tilde p_1,\dots,\tilde p_6$ is in normal form:
\begin{eqnarray*}
\tilde p_1 = u_1^5  - u_1^3 u_2; \quad
\tilde p_2 = u_2^7 - u_1^8  u_2;\\
\tilde p_3 = u_3^2 - u_1^6 u_3^3;\quad
\tilde p_4 = u_4^3 - u_1^2 u_2 u_5^2 u_6;\\
\tilde p_5 = u_5^2 - u_4^2 u_5 ;\quad
\tilde p_6 = u_6^3 - u_1^3 u_4 u_5 u_6.
\end{eqnarray*}
Thus, we conclude that $p_1,\dots,p_8$ defines a complete intersection.
We will compute the numerical invariants of this system in 
Example~\ref{excont}.
\end{example}

\section{Computing the number of solutions}
\label{sec:irreducible}

We recall that if $\pn$ is a gci then we denote by $d$ (respectively $D$)
the number of points in $\VV_c \cap \bk^n$ counted
with multiplicity (respectively without multiplicity), for a generic choice of non-zero
coefficients.  Similarly, recall that for
any index set $L \subset \{1,\dots,n\}$ we denote by
  $\mu_L$ the number of points in $\VV_c \cap \bk^n_L$ counted
with multiplicity,
where $\bk^n_L$ is the set of points
in affine space whose coordinate $x_\ell=0$ precisely when $\ell \in L$. 
In particular, $\mu=\mu_{[n]}$ denotes the multiplicity
at the origin.

If  $\pn$ is a gci but $0\not\in \VV(\JJ)$, then it follows from
Lemma~\ref{preinv} and Proposition~\ref{invertible} that the
invariants $d$ and $D$ of $\pn$ are obtained from those of the
derived system by multiplying times $|\det B_2|$.  
We will  assume from now on that no variable is invertible
modulo $\JJ$, i.e., that 
$0\in \VV(\JJ)$.

\medskip

We begin this section by showing that  it is enough to compute
the desired numerical invariants $d, D, \mu_L$, for 
ideals in normal form.    We then show that if the system is irreducible, 
in a sense
made precise below, then the only zero outside the torus is the
origin and its multiplicity may be easily computed from the 
exponents of the system.  Finally, we consider the general case and
show how the various dimensions depend on the combinatorics of
the irreducible components.  

\subsection{Multiplicities and parametric reduction.}

  Suppose $\pc$ is as
in (\ref{thesystem}) with $p_n =  x_n^{\ell} - c_n x_{n-1}^{m} $,
$\ell, m > 0$.  Let
$q = {\rm gcd}(\ell,m)$ and
$$p'_n  \ =\ x_n^{\ell'} - c'_n x_{n-1}^{m'}\,.$$
We will denote by $d',D',\mu'_L$ the corresponding invariants
for $p_1,\dots,p_{n-1},p'_n$.
  
We show, first of all, that by keeping track 
of
$q$ we may assume without loss of generality that $m$ and $\ell$
are coprime.

\begin{lemma}\label{multred}
With notation as above, set $m' = m/q$, $\ell'=\ell/q$,
$p'_n =  x_n^{\ell'} - c'_n x_{n-1}^{m'} $ and let $B$ and $B'$ be
the corresponding matrices.
\begin{enumerate}
\item $|\det B| \ =\ q\cdot|\det B'|$.
\item $\pn$ is a gci if and only if
$p_1,\dots,p_{n-1},p'_n$ is a gci.
\item For any index set $L \subset \{1,\dots,n\}$,
$\mu_L \ =\ q\cdot\mu'_L$.
\item $d = q\cdot d'$ and $D = q\cdot D'$.
\end{enumerate}
\end{lemma}

\begin{proof} The first assertion is trival while the second one
follows from Theorem~\ref{fscriterion}.  In order to
prove assertion $3$, let $(c_1,\dots,c_n)\in (k^*)^n$ be such that 
$\JJ_c$ is a complete intersection and decompose
\begin{equation} \label{eq:fact}
p_n =x_n^{\ell} - c_n x_{n-1}^{m}  = \prod_{\xi\in W_q} ( 
x_n^{\ell'}- \xi \, x_{n-1}^{m'})\,,
\end{equation}
where $W_q$ denotes the $q$-th roots of $c_n$. For any $\lambda \in \VV_c$, 
we have
$$\dim_k \left(R_\lambda/(\JJ_c)_\lambda\right) \ =\ \sum_{\xi\in W_q} \dim_k 
\left(R_\lambda/(\JJ_\xi)_\lambda\right)\,,$$
where $\JJ_\xi:= \langle p_1(c;x), \dots, p_{n-1}(c;x),x_n^{\ell'}
-\xi\,x_{n-1}^{m'} \rangle$.  Therefore,
$$\sum_{\lambda \in \VV_c \cap \bk^n_L}\dim_k 
\left(R_\lambda/(\JJ_c)_\lambda\right) \ =\ \sum_{\xi\in W_q} \
\sum_{\lambda \in \VV(\JJ_\xi) \cap \bk^n_L}\dim_k 
\left(R_\lambda/(\JJ_\xi)_\lambda\right)\,,$$
By a scalar change of variables it follows that 
$$\sum_{\lambda \in \VV(\JJ_\xi) \cap \bk^n_L}\dim_k 
\left(R_\lambda/(\JJ_\xi)_\lambda\right)\,,$$ 
is independent of $\xi\in W_q$ and, since it agrees with $\mu'_L$, we
obtain that
$$\mu_L\ =\ q\, \mu'_L$$
as claimed.  The last assertion follows directly from the previous one
and the factorization (\ref{eq:fact}).
\end{proof}

We next show that multiplicities are not altered
under parametric reduction.  If  the binomial system $\pn$
is a gci, and
$p_n = 
x_n^{\ell} - c_n x_{n-1}^{m} $,
$\ell, m > 0$ coprime, let
$\tildep$ be the binomial system obtained through
parametric reduction.  We will denote by
$\tilde d, \tilde D$ and $\mu_{\tilde L}$ the
corresponding invariants.

Given $L \subset [n]$ we denote by
$\tilde L := L \cap [n-1]$.  Conversely, given
$\tilde L \subset [n-1]$ set $L = \tilde L$ if $n-1 \not\in \tilde L$ and 
$L = \tilde L \cup \{n\}$ otherwise.  Note that
if $L\subset [n]$ is such that $\mu_L\not=0$ then
either $L\subset [n-2]$ or both $n-1,n \in L$.  Hence,
the correspondence $L \mapsto \tilde L$ 
establishes a bijection between index sets $L\subset [n]$
such that $\mu_L \not = 0$ and subsets 
$\tilde L\subset [n-1]$ such that $\mu_{\tilde L}\not = 0$.

\begin{theorem}\label{parred}
Suppose that $\pn$ is a gci and
 $p_n = 
x_n^{\ell} - c_n x_{n-1}^{m} $, with
$\ell, m $ coprime positive integers.
Let $\tildep$ be the binomial system obtained through
parametric reduction. Then $D = \tilde D$ and, for any $L\subset [n]$,
\begin{equation}\label{dim}
\mu_L = \mu_{\tilde L}.
%\dim_k R/\JJ \ =\ \dim_k \tilde R /\tilde \JJ.
\end{equation}
Consequently, $d=\tilde d$ as well.
\end{theorem}

\begin{proof}  
Let $c=(c_1,\dots,c_n)\in (k^*)^n$ be such that 
$\JJ_c$ is a complete intersection.  We may assume without
loss of generality that $c_n=1$.  Let $\tilde c = (c_1,\dots,c_{n-1})$
and denote by $\tilde \JJ_{\tilde c}$ the ideal generated by
$\tilde p_1(\tilde c;u), \dots, \tilde p_{n-1}(\tilde c;u)$ in the
ring $k[u]$.
Given any
$\tilde\lambda = (\lambda_1 , \dots,  \lambda_{n-1})\in
\VV(\tilde \JJ_{\tilde c}) \subset \bk^{n-1}$, let us denote by
$\lambda$ the point $ (\lambda_1 , \dots, \lambda_{n-2},
\lambda_{n-1}^\ell,
\lambda_{n-1}^m)\in
\VV(\JJ_c) \subset \bk^{n}$.  This assignment $\tilde\lambda \mapsto \lambda$
defines a bijection between $\VV(\tilde \JJ_{\tilde c})$ and $\VV(\JJ_c) $
since $\ell, m$ are coprime, and so $D = \tilde D$.  To show that 
$d = \tilde d$ it suffices to prove that
at the level of local rings
\begin{equation}\label{dimlocal}
\dim_{\bk} (R\otimes_k \bk)_\lambda/(\JJ_c)_\lambda \ =\ \dim_{\bk} (
\tilde R\otimes_k \bk)_{\tilde \lambda} 
/(\tilde \JJ_{\tilde c})_{\tilde \lambda}.
\end{equation}

We will denote by
$A_1$ the localization of $\bk[u_1,\dots,u_{n-1}]$ at
$\tilde\lambda$ and by $A_2$ the localization  of
$\bk[u_1,\dots,u_{n-2},u_{n-1}^\ell,u_{n-1}^m]$ at $\tilde\lambda$.  Let
$(\hat \JJ_{\tilde c})_{\tilde\lambda}$ be the ideal generated by 
$\tilde p_1(\tilde c;u),\dots,\tilde p_{n-1}(\tilde c;u)$ in $A_2$ so that
$(\tilde \JJ_{\tilde c})_{\tilde\lambda} = A_1\cdot (\hat \JJ_{\tilde c})_{\tilde\lambda}$.  
Again, since $m$ and $\ell$ are coprime it is clear that
$$
\dim_{\bk} (R\otimes_k \bk)_\lambda/(\JJ_c)_\lambda \ =\ \dim_{\bk} A_2/(\hat \JJ_{\tilde c})_{\tilde\lambda}\,.
$$
Thus, the result will follow if we show that
\begin{equation}\label{eq1}
\dim_{\bk} A_1/(\tilde \JJ_{\tilde c})_{\tilde\lambda} \ =\
  \dim_{\bk} A_2/(\hat \JJ_{\tilde c})_{\tilde\lambda}
  \end{equation}
The following proof of (\ref{eq1}) was suggested to us by Mircea 
Mustata.

We recall from \cite[\S14]{matsumura} the following notion
of {\em multiplicity}:  Let $(R,\mm)$ be a $d$-dimensional
Noetherian local ring, $M$ a finite $R$-module and
$\qq$ an $\mm$-primary ideal.  The multiplicity of
$M$ with respect to $\qq$ equals
\begin{equation}\label{multiplicity}
\ee(\qq,M)\ =\ \lim_{m\to\infty}\,
\frac{d!}{m^d} \,{\rm length}(M/\qq^{m+1}M)
\end{equation}

Since both $A_1$ and $A_2$ are Cohen-Macaulay rings of dimension
$n-1$ and $$\tilde p_1(\tilde c;u),\dots,\tilde p_{n-1}(\tilde c;u)$$ define a regular sequence in $A_2$, hence in
$A_1$ as well, it follows from \cite[Theorem~14.11]{matsumura}
that
\begin{equation*}\label{m1}
\dim_{\bk} A_2/ (\hat \JJ_{\tilde c})_{\tilde\lambda} \ =\ \ee( (\hat \JJ_{\tilde c})_{\tilde\lambda}, A_2)\ \
\hbox{and}\ \
\dim_{\bk} A_1/(\tilde \JJ_{\tilde c})_{\tilde\lambda} \ =\ \ee((\tilde \JJ_{\tilde c})_{\tilde\lambda}, A_1).
\end{equation*}
On the other hand, $A_1$ may be considered as a $A_2$-module
and it is clear from (\ref{multiplicity}) that
$$ \ee((\tilde \JJ_{\tilde c})_{\tilde\lambda}, A_1)\ =\ \ee( (\hat \JJ_{\tilde c})_{\tilde\lambda}, A_1)\,.$$
Finally, \cite[Theorem~14.8]{matsumura}
gives that
$$\ee( (\hat \JJ_{\tilde c})_{\tilde\lambda}, A_1) = {\rm rank}_{A_2} A_1\cdot \ee( (\hat \JJ_{\tilde c})_{\tilde\lambda},A_2) =
  \ee( (\hat \JJ_{\tilde c})_{\tilde\lambda}, A_2),$$
  since the assumption that $m$ and $\ell$ are coprime
  implies that  the two domains $A_1, A_2$ have
  the same fraction field and so ${\rm rank}_{A_2}A_1 =1$.  
This proves (\ref{eq1}).
\end{proof}

\subsection{Irreducible Systems}

\begin{definition} \label{def:irreducible}
A binomial system $\pn$ is said to be
irreducible if it is in normal form and it is not possible
to reorder it so as
to find a proper index subset $I\subset [n]$ such that for every $i\in I$ the
binomial $p_i$ depends only on the variables $x_j, j\in I$.
\end{definition}

Recalling that a system in normal form is a gci and that
$0\in \VV(\JJ)$, we easily have:

\begin{lemma} \label{lemma:irredzero}
Let $\p$ be an irreducible system as in (\ref{thesystem}) and
let $c\in (k^*)^n$ be such that $\JJ_c$ is a
complete intersection. 
Then if $a\in \VV_c$,  either
$a=0$ or $a\in (\bk^*)^n$.
\end{lemma}

\begin{proof} Given $a \in \VV_c$, let
$I = \{i\in [n] : a_i \not= 0\}$.
If $i\in I$ then, since $\p$ is in normal form,
$$ p_i (c;x)= x_i^{r_i} -c_i x^{\beta_i}\,;\ r_i>0,\ \beta_i\not=0,$$ 
and, since $a_i \not=0$, we
must have ${\rm supp}(\beta_i) \subset I $ for all $i\in I$.
This contradicts the irreducibility of $\p$ unless
$I=[n]$ or $I=\emptyset$.
\end{proof}

The following theorem identifies $d$ and $\mu$ for irreducible systems.
Recall that $\delta = |\det B|$ is the cardinality
of $\VV_c \cap (\bk^*)^n$.
Our arguments are built on
the proof of a result of Vinberg (cf. \cite[Theorem 4.3]{kac}).

\begin{theorem}\label{vinberg}
Given an irreducible system 
\begin{equation*}
 p_i(c;x) \ =\ x_i^{r_i}- c_i x^{\beta_i}\,, \ 
i=1,\dots,n,
\end{equation*}
where $ r_i>0,\  \beta_i \in \N^n,\ \beta_i \not= 0$, then:
\begin{itemize}
\item If all principal minors of $B$ are positive
  $$d = r_1\cdots r_n\ ;\quad \mu = d - \delta.$$
Such a system will be called a {\em global} irreducible system.

\item Otherwise, $\mu = r_1\cdots r_n$ and $d = \mu + \delta $.
In this case we say that the system is {\em local}.
\end{itemize}
\end{theorem}

\begin{proof}
Let us fix throughout coefficients $c\in (k^*)^n$ such
that $\JJ_c$ is a complete intersection.
Since the system is in normal form, the entries of 
$B$ are $b_{ii} = r_i - (\beta_i)_i$ and
$b_{ij} = - (\beta_i)_j, i\not=j$. Hence, its off-diagonal
terms are non-positive.  Moreover, the irreducibility of
the system implies that $B$ is indecomposable in the
sense of \cite {kac}.  In fact, the irreducibility of
the system implies a stronger condition, namely
\cite[Lemma~4.3]{kac}: Suppose $u\in \R^n$ is
a vector with non-negative entries and that $B\cdot u \geq 0$ in
the sense that all its entries are non-negative as well.  Then
either $u=0$, or $u>0$, i.e., all its entries are strictly positive.
Indeed, let $I= \{i\in [n] : u_i = 0$, then for any $i\in I$,
$(B\cdot u)_i \leq 0$ and equality occurs if and only
if $b_{ij}=0$ for all $j\not\in I$.  Hence, by irreducibility
we must have $I = [n]$ or $I = \emptyset$.

Given that \cite[Lemma~4.3]{kac} holds in our case, we can
apply Theorem~4.3 in \cite{kac} and conclude that three cases
are possible:

\begin{itemize}
\item There exists  $w\in \Q^n$ all of whose entries
are positive such that $B\cdot w >0$.

\item There exists $w\in \Q^n$,  all of whose entries
are positive such that
$B\cdot w <0$.

\item ${\rm rank}(B) = n-1$ and there exists 
$w\in \Q^n$ all of whose entries
are positive such that $B\cdot w =0$.
\end{itemize}
According to \cite[Theorem 2.3]{bp}, the first condition
is equivalent to the statement that all principal minors
of $B$ are positive which implies, in particular,
that all the diagonal entries of $B$ are strictly
positive.  These are the so-called $M$-matrices
of \cite{bp}.  
Moreover, if we consider a term order in $\kn$ that
refines the weight order defined by $w$, the term $x_i^{r_i}$
will be the leading term in $p_i(c;x)$, and hence $\pc$ is a Gr\"obner
basis.  It then follows \cite[\S 5.3, Proposition~4]{clo1} that  $d = r_1\cdots r_n$ and, by 
Lemma~\ref{lemma:irredzero}, $\mu = d - \delta$.

In the second case we can similarly
define a local order (cf. \cite{singular}) for which the leading term of $p_i(c;x)$ is $x_i^{r_i}$.
Hence $\pc$ is a standard basis in the local quotient ring at the origin 
and, consequently,  $\mu = r_1\cdots r_n$ and $d = \mu + \delta$.
We note that this is valid even if $\det B=0$ since,
in that case, $\JJ_c$ is a complete intersection if and only
if  $\VV_c = \{0\}$.

In the third case, the binomials $p_i(c;x)$ are weighted homogeneous
relative to the weight $w$ and therefore $\mu = r_1\cdots r_n$ and
$d = \mu + \delta$ since, again, $\VV_c$ consists of only the origin.
Thus this case behaves as the previous one and we will also refer
to it as a local case.
\end{proof}

\begin{remark}\label{local1}
We note that 
if $n=1$, the system $p = x^\alpha - c x^\beta, \alpha\not=\beta$, will be local
if $\alpha < \beta$ and global if $\alpha > \beta$.
\end{remark}

\subsection{The General Case}

We consider now general gci systems in normal form.
Throughout this subsection we will, again, fix 
coefficients $c\in (k^*)^n$ so that $\JJ_c$ is
a complete intersection.  For economy of notation
we will denote simply by $p_i$ the corresponding
binomials in $\kn$.
If the system $\p$ is not irreducible, 
then, as Lemma~\ref{tarjan} shows, it is possible
to choose an increasing sequence
\begin{equation}\label{nus}
0= \nu_0 < \nu_1 < \cdots < \nu_s = n
\end{equation}
so that if $I_a = \{\nu_{a-1}+1,\dots,\nu_a \}$, then the following
holds:
\begin{itemize}
\item For $i\in I_a$, $p_i \in k[x_j; j\in I_1\cup \cdots \cup I_a]$.
\item The system
$
\hat p_i := p_i(1,\dots,1,x_{\nu_{a-1}+1},\dots,x_{\nu_a })
$,
$i\in I_a$, is irreducible.
\end{itemize}

\begin{definition} \label{def:triangular}
A system of this form will be said to be in {\em triangular} form
relative to the blocks $I_1,\dots, I_s$.  Given a reducible system in triangular form, 
we will refer to the system
$\{\hat p_i, i\in I_a\}$ as the restriction of $\p$ to
$I_a$ and denote it, for short, by $\hat p^a$.
\end{definition}

\begin{lemma}\label{tarjan}
Any system of $n$ binomials $\p$ in normal form (\ref{def:nf})
can be put in triangular form  in time $O(n^2)$.
\end{lemma}

\begin{proof}
Consider the ocurrence matrix $N$: this is a $0$-$1$ matrix
with $n_{ij} \not= 0 $ if and only if  $i \not= j$ and $p_i$ depends on $x_j$
 (i.e., if $p_i = x_i^{r_i} -c_i x^{\beta_i}$ with $\beta_{ij}
\not= 0$).
This is a standard construction, first used by Steward \cite{steward},
for the analysis of the structure of large systems of equations.
Note that, because the system is in normal form,  putting $\p$ in triangular form
corresponds precisely to finding a permutation matrix $P$ such that
${}^tP N P$ is block lower triangular, with the irreducible subsystems of $\p$
corresponding to  the irreducible diagonal square 
blocks along the diagonal of $\,{}^tP N P$.

Tarjan's algorithm \cite{tarjan} to search for the strongly connected
components of the directed graph associated to $N$ provides an efficient method 
for finding such permutation matrix $P$ \cite{dr, bahia}; it runs in time
linear in the number of vertices plus the number of edges of the graph. 
\end{proof}

Given a system in normal form and triangular relative to $I_1,\dots, 
I_s$, let
 $\delta_a = |\det  B_a|$ , where  $B_a$ is the matrix
associated with the
system $\hat p^a$ and
$$\rho_a = \prod_{j\in I_a} r_j.$$
 We also denote by $\mu_a$ the multiplicity of $\hat p^a$
 at $0$ and by
$d_a$ the total number of solutions of $\hat p^a$
counted with multiplicity. 

For a triangular system $\p$, its associated matrix is
block lower-triangular:
\begin{equation}\label{triangular}
B = \left(
\begin{array}{cccc}
B_1 & 0 & \dots & 0\\
C_{21} & B_2 & \dots &0\\
\vdots & \vdots & \ddots &\vdots\\
C_{s1}&C_{s2}&\dots & B_s
\end{array}
\right).
\end{equation}

The number of solutions of the system $p_1, \dots, p_n$ and the patterns of possible
zero
coordinates of the solutions are best described in
terms of  the directed acyclic graph $G$ with $s$ vertices labeled $\{1, \dots, s\}$
and an arrow from node $a$ to node $b$ if and only if the rectangular
submatrix $C_{ba}$ is not identically zero. 
We recall that a vertex  is called
a {\em source \/} if it is not the head of any arrow. The subset of
sources of the vertex set of a subgraph $H$ of $G$ will be denoted by $S(H)$.

\begin{remark} 
We can think of $G$ as a weighted graph, where each vertex $a \in [s]$ comes
 with the weights  $\delta_a, \rho_a, \mu_a$ (or $\delta_a, d_a, \mu_a$). 
Equivalently, we can think  that the information at each node
is coded by the weights $\delta_a, \rho_a$ plus an additional label
{\em global\/}  or {\em local}
according to where $B_a$ is global or local, which prescribes the relation among 
$\delta_a, \rho_a$ and $\mu_a$ (or $\delta_a, d_a$ and $\mu_a$). 
\end{remark}

\begin{theorem}\label{th:localmult}
The multiplicity  $\mu$ of $\JJ_c$ at the origin
equals
\begin{equation} \label{eq:multat0}
\mu \, = \, \left( \prod_{a \in G \backslash S(G)} \rho_a \right) \, \left( \prod_{b \in S(G)} \mu_b \right).
\end{equation}
\end{theorem}

\begin{proof}
We will prove formula~\ref{eq:multat0} by induction in the
number $s$ of blocks.  If $s=1$, the system is irreducible and
$\{1\}\in S(G)$ so the formula holds. Consider $s > 1$ and assume
that the result is true for systems with $s-1$ blocks.
Let $B$ be as in (\ref{triangular}), set $n':=\nu_{s-1}$, where
$\nu_{s-1}$ is as in (\ref{nus}), and consider the ideal
$\JJ'_c:= \, \langle p_1, \dots, p_{n'} \rangle\, , $
in the polynomial ring in the first $n'$ variables.  Clearly, $p_1,
\dots, p_{n'}$  is in normal and triangular form.
Let $G'$ be the corresponding graph; it is obtained by erasing from $G
$ the vertex $s$ and
all edges ending at $s$. By inductive hypothesis, we have
that the multiplicity $\mu'$ of $\JJ'_c$ at $0'$ equals
\begin{equation}\label{eq:multprima}
\mu' \, = \, \left( \prod_{a \in G' \setminus S(G')} \rho_a \right)
\, \left( \prod_{b \in S(G')} \mu_b \right).
\end{equation}
The matrix $B$ has the form
\begin{equation}\label{Bprime}
B = \left(
\begin{array}{ccccc}
{}& {} & {} & \vrule & {}\\
{} & B' & {} & \vrule &0\\
{}& {} & {} & \vrule & {}\\
\hline
{}&C&{} &\vline & B_s
\end{array}
\right).
\end{equation}
If the rectangular matrix $C$ is identically zero, then the last $n-
n'$ polynomials
depend only on the last $n-n'$ variables, and we  have that
$$ \mu \, = \mu' \, \cdot \, \mu_s,$$
as wanted, since in this  case $S(G) = S(G') \cup \{s\}$.

On the other hand, if $C$ is not zero, it is possible to find a
positive weight vector
$w$ such that the initial monomial $in_{-w}(p_j) = x_j^{r_j}$, for
all $ n' < j \leq n$.  Indeed,  set  $J_0 = [n'] \,$  and,
for $l\geq 1$ define
$$J_{\ell} \,:=\,\left\{k \in [n] \backslash \left(\bigcup_{a=0}^{\ell-1}
J_a\right)\, : \,J_{\ell-1} \cap {\rm supp}(\beta_k) \not =\emptyset
\right\}.$$
Note that $C\not=0$ implies that $J_1$ is non empty.  
Also, the assumption that 
$B_s$ is irreducible guarantees that there exists $L \leq n- n'$ such
that  $[n]
\backslash [n'] = \bigcup_{1 \leq \ell \leq L} J_{\ell}$.  
Now, choose $w_k =1 $ for $k\in J_L$.  Then assuming that the
weights for the variables $k\in J_a$, $\ell \leq a \leq L-1 $, have been chosen
so that $in_{-w}(p_j) = x_j^{r_j}$ for all $j\in J_b$, $b\geq \ell+1$, 
we may choose positive weights $w_k$ for $k\in J_{\ell-1}$ that are
sufficiently large so that  
$in_{-w}(p_j) = x_j^{r_j}$ for all $j \in J_\ell$ as well.

Consider now any local order $\prec$ in $k[x_1,\dots, x_n]$ refining
the weight $-w$.
Let  $\{q_1,\dots,q_t\}$ be a standard basis for
the ideal $\JJ'_c$ with respect to the local order induced by $\prec$
in $k[x_1,\dots,x_{n'}]$. Then,
$\{q_1,\dots,q_t,p_{n'+1},\dots,p_n\}$
is a standard basis for $\JJ_c$ relative to $\prec$ since,
for every $i=1,\dots,t$, the
leading monomials of the polynomial $q_i$ is
coprime with those of the $p_j$, $n' < j \leq n$,
and, therefore, the weak normal
form of the corresponding $S$-polynomial is $0$ \cite{singular}.
The corresponding
initial ideal $L_{\prec}(\JJ_c)$ will be generated by some monomials  in
the first $n'$ variables (generating the initial ideal $L_{\prec'}
(\JJ'_c)$)
and the pure powers $x_j^{r_j}$ for all $j > n'$.
Therefore, the multiplicity $\mu$ of $\JJ_c$ at $0$ equals:
\begin{eqnarray*}
\dim_{\bk}
\left( {\bk}[x_1 \dots, x_n] / \JJ_c \right)_{0} =
\dim_{\bk}
\left( {\bk}[x_1 \dots, x_n] / L_{\prec}(\JJ_c) \right)_{0} = \\
\dim_{\bk}
\left( {\bk}[x_1 \dots, x_{n'}] / L_{\prec'}(\JJ'_c) \right)_{0} \,
\cdot \,
\dim_{\bk}
   \left( {\bk}[x_{n'+1} \dots, x_{n}] / \langle x_{n'+1}^{r_{n'+1} }
\dots x_n^{r_n} \rangle \right)_{0} =\\
\dim_{\bk}
\left( {\bk}[x_1 \dots, x_{n'}] / \JJ'_c  \right)_{0} \, \cdot \,
\rho_s.
\end{eqnarray*}
In this case $s\not\in S(G)$, and so $S(G) = S(G')$. Since
the dimension of the local quotient by $\JJ'_c$ at the origin equals
(\ref{eq:multprima}), we get that
$$\mu = \mu'  \cdot \rho_s = \, \left( \prod_{a \in  \setminus S(G)}
\rho_a \right) \, \left( \prod_{b \in S(G)} \mu_b \right),$$
as wanted.
\end{proof}

\begin{remark} 
Using Theorem~\ref{vinberg}  we can translate (\ref{eq:multat0}) as
\begin{equation} \label{eq:multat0v2}
\mu \, =  \, \left( \prod_{ a \in G_1 } d_a \right) \, \left( \prod_{b \in G_2} \mu_b \right),
\end{equation}
where $G_1$ is the set of nodes of $G$ corresponding to the
global, non-sources of $G$ and $G_2$ is its complement.
\end{remark}

\medskip

We will also need the following terminology.

\begin{definition} \label{def:full}
A vertex $b$ of (the directed acyclic graph) $G$ is said to be a {\em descendant\/}
 (respectively, a {\em direct descendant}) of the vertex $a$ if there is
a directed path (respectively, a directed edge) from $a$ to $b$. 
  A (directed) subgraph $H$ of $G$
is said to be {\em full\/} if, for any of its vertices $j$, all its descendants and 
all the directed paths starting from $j$
 also belong to $H$. 
The collection of full subgraphs of $G$ will be denoted by 
$\FF(G)$.
\end{definition}

The empty subgraph is full and even if $G$ is connected, a full subgraph $H$ may be disconnected.
Note also that a full subgraph is completely determined by its sources.

The following result refines the description given in Remark~\ref{setsL} of subsets $L\subset [n]$
with $\mu_L \not =0$ .

\begin{proposition} \label{pattern}
Let $\p$ be a binomial complete intersection in normal and triangular form and $L\subset [n]$.
Then $\mu_L = 0$ unless 
there exists a full subgraph $H$ of $G$ such that
\begin{equation}\label{detnonzero} 
\prod_{a \notin H} \delta_a \not = 0
\end{equation}
and $L$ coincides with the union of all the indices belonging to
blocks that are vertices of $H$. 
\end{proposition}

\begin{proof}
With the above notations, let $\lambda = (\lambda_1, \dots, \lambda_n) \in \VV(\JJ_c)$ and 
$L = L(\lambda) = \{ i \in [n] \, : \, \lambda_i=0\}$.
Set $H = \{a\in G : I_a \cap L \not= \emptyset\}$.
If $a\in H$ then we may argue as in 
 Lemma~\ref{lemma:irredzero} to conclude that $I_a \subset L$.
 Suppose now that $a\in H$ and that $(a,b)$ is an edge
 in $G$.  Since $C_{ba} \not= 0$, there exists $i\in I_a$ and
 $j\in I_b$ such that $i\in {\rm supp}(\beta_j)$ and, consequently,
 $\lambda_j = 0$, i.e., $j\in I_b \cap L$, and $b\in H$.
This shows that $H$ is a full subgraph of $H$. 
The need for condition (\ref{detnonzero}) was already noted in Remark~\ref{setsL}.
\end{proof}

With notation as in Prop.~\ref{pattern}, given a full subgraph
$H\subset  G$, we will denote by $L(H)$ the set of indices
belonging to blocks associated with vertices of $H$.

\begin{proposition} \label{prop:muL}
Given a full subgraph $H$ of $G$, the number $D_{L(H)}$ of points in 
$\VV(\JJ_c) \cap \bk^n_{L(H)}$ counted without multiplicity
equals
\begin{equation} \label{eq:DL}
D_{L(H)} \, = \, \left(\prod_{a \notin H} \delta_a \right)
\end{equation}
while the number $\mu_{L(H)}$ of points in $\VV(\JJ_c) \cap \bk^n_{L(H)}$
counted
 with multiplicity equals 
\begin{equation}\label{muL}
\mu_{L(H)} \, = \, \left(\prod_{a \notin H} \delta_a \right) \, \left( \prod_{b \in H \setminus S(H)} 
\rho_b \right) \, \left( \prod_{e \in S(H)} \mu_e \right).
\end{equation}
\end{proposition}

\begin{proof}
The first assertion follows easily from Proposition~\ref{pattern}.
In order to prove (\ref{muL}), let
$\lambda \in \VV(\JJ_c) \cap \bk^n_{L(H)}$,  write $\lambda = ( \lambda^{(1)}, \dots, \lambda^{(s)})$ with
$\lambda^{(a)} \in  {\bk}^{|I_a|}$ for all $ a \in [s]$. 
Since $H$ is a full subgraph, there are no edges starting
at a node in $H$ and ending at a node outside of $H$; i.e., $C_{ba}=0$ for all $a \in H$ and $b \notin H$.
Therefore, it is possible to relabel the variables and the
binomials $\pn$ so that the system remains in normal form
and satisfies that  $a < b$ for
all $a \notin H$ and $b \in H$. Thus, we may assume
without loss of generality that $H = \{t+1, \dots, s\}$ and
therefore
$\lambda = (\lambda^{(1)}, \dots, \lambda^{(t)}, 0, \dots, 0)$ with
$\lambda^{(a)} \in  ({\bk}^*)^{|I_a|}$ for $a = 1, \dots, t$.  
Equivalently, $$\lambda = (\lambda',0) \in ({\bk}^*)^{n'}\times({\bk})^{n-n'}\ ;\ n' := \nu_t\,.$$
We let $x'$ stand for the first $n'$ variables $x_1,\dots,x_{n'}$ and
$x''$ for the remaining $n-n'$ variables.  Then 
$$\JJ'_c \ :=\ \langle p_1,\dots,p_{n'}\rangle \subset k[x']$$ and $\lambda'$ 
is a simple zero of $\JJ'_c$.  Hence $p_1,\dots,p_{n'}$ define
the maximal ideal in the local ring $\left({\bk}[x']\right)_{\lambda'}$.
We then have:
\begin{eqnarray*}
 \mu_\lambda &:= &\dim_{\bk} \left( {\bk}[x] /  \JJ_c\right)_{\lambda} \\
& = & \dim_{\bk} \left( {\bk}[x] /  \langle  x_1 - \lambda_1,\dots, x_{n'} - \lambda_{n'}, p_{n'+1}, \dots, p_n
\rangle\right)_{\lambda} \\
& = & \dim_{\bk} \left( {\bk}[x''] /  \langle p_{n'+1}(\lambda',x''), \dots, p_n(\lambda',x'')
\rangle\right)_{0}\\
&= & \dim_{\bk} \left( {\bk}[x''] /  \langle p_{n'+1}(1, \dots, 1,x''), \dots, p_n(1, \dots, 1,x'')
\rangle\right)_{0}.
\end{eqnarray*}
So,  $\mu_\lambda$  equals the multiplicity at the origin
$0 \in {\bk}^{n-n'}$ of the system $\{\hat p_{n'+1}, \dots, \hat p_n \}$.
Formula (\ref{muL}) now follows from Theorem~\ref{th:localmult}, 
and the fact that the system $p_1, \dots, p_{n'}$ has
 $\delta_1 \cdots \delta_t$ simple solutions in $({\bk}^*)^{n'}$.
\end{proof}

The following explicit formulas for $d$ and $D$ follow
by adding (\ref{eq:DL}) and (\ref{muL}) over all full subgraphs
of $G$.

\begin{theorem} \label{th:dformula}
Suppose that
$\pn$ are in normal, triangular form. For generic parameters
$c \in (k^*)^n$, the total
number of solutions of the system $p_1(c;x) = \cdots =
p_n(c;x) =0$, counted without multiplicity, equals
\begin{equation}\label{eq:Dformula}
D \, = \, \sum_{H \in \FF(G)} \left(\prod_{a \notin H} \delta_a \right),
\end{equation}
and the total number of solutions counted with multiplicity equals
\begin{equation}\label{eq:dformula}
d \, = \, \sum_{H \in \FF(G)} \left(\prod_{a \notin H} \delta_a \right) 
\, \left( \prod_{b \in H \setminus S(H)} 
\rho_b \right) \, \left( \prod_{e \in S(H)} \mu_e \right).
\end{equation}
\end{theorem}

We end this section with a recursive formula to compute $d$.
In order to state the following proposition we define, for $1\leq r \leq
s$, the
binomial system $q^{(r)}$:
$$ p_i(1,\dots,1,x_{\nu_{r-1}+1},\dots,x_n) \ ,  \quad i \in I_r \cup \dots \cup I_s.$$
Note that the matrix associated with $q^{(r)}$ is:
\begin{equation}\label{ktriangular}
B^{(r)} = \left(
\begin{array}{cccc}
B_r & 0 & \dots & 0\\
C_{(r+1)r} & B_{r+1} & \dots &0\\
\vdots & \vdots & \ddots &\vdots\\
C_{sr}&C_{s(r+1)}&\dots & B_s
\end{array}
\right)
\end{equation}
Clearly if $\pn$ is in normal, triangular form, so is $q^{(r)}$.
We denote by $F_r$ the number of solutions in $\bk^{ n - \nu_{r-1}}$, 
counted with multiplicity, of the system $q^{(r)}$.

\begin{proposition}\label{prop:inductive}
$F_r$ is
a polynomial function of $\{\delta_a, \mu_a, \rho_a\,;\ a= r,\dots,s\}$. 
It may be computed recursively as:
$$
F_s = d_s = \delta_s + \mu_s$$
\begin{equation} \label{eq:Fk}
 F_r = \delta_r \cdot F_{r+1} + \mu_r \cdot
F_{r+1}|_{\delta_b=0,\mu_b=
\rho_b},
\end{equation}
where $b$ runs over all indices in $\{r+1,\dots,s\}$ such that $C_{b r} \not=0$.
\end{proposition}

\begin{proof}
We may assume without loss of generality that  $r =1 < s$.
Let $G$ be the graph of $B$ and
 $G^{(2)}$  the subgraph of $G$ associated to the 
submatrix $B^{(2)}$ defined by (\ref{ktriangular}).

Any full subgraph $H\in \FF(G^{(2)})$ 
may be thought of as a full subgraph in $G$.  We denote
by $\FF' \subset \FF(G)$ the collection of such subgraphs.
Clearly $\FF'$ consists of all full subgraphs of $G$ not
containing the vertex $1$.  Let $\FF''$ denote the complement
of $\FF'$ in $\FF(G)$.  Removing the vertex $1$ from a subgraph
$H\in \FF''$ defines a full subgraph $H^{(2)}$ of $G^{(2)}$
with the property that no direct descendant of $1$ in $G$ may be in
$G^{(2)}\setminus H^{(2)}$.  Let us denote by $\FF''(G^{(2)})$
the collection of such full subgraphs of $G^{(2)}$.
We can write 
\begin{equation}\label{sum}
F_1 \ =\ \sum_{H\in \FF'} \mu_{L(H)} +  \sum_{H\in \FF''} \mu_{L(H)}.
\end{equation}
Since, for $H\in \FF'$, $1\not\in H$, in view of 
(\ref{muL}), the first sum may be computed
as:
\begin{equation}\label{primera}
\sum_{H\in \FF'} \mu_{L(H)}\ =\  \delta_1 \sum_{H\in \FF(G^{(2)})} \mu_{L(H)} \ =\  \delta_1 F_2,
\end{equation}
since  $S(H)$ is the same whether we view 
$H$ as a subgraph of $G$ or of $G^{(2)}$.

Thus, in order to complete the proof we need to show that the
second sum in (\ref{sum}) equals 
$$\mu_1 \cdot
F_{2}|_{\delta_b=0,\mu_b=
\rho_b},$$
where $b$ runs over all vertices in $G^{(2)}$ that are direct
descendants of $1$ in $G$.  We note first of all, that  setting 
$\delta_b=0$ for all direct descendants $b$ of $1$ 
has the effect of restricting the sum in (\ref{eq:dformula}) to  $\FF''(G^{(2)})$.  
Moreover, given $H\in \FF''$, let
$H^{(2)}$ denote the full subgraph of $G^{(2)}$ obtained by removing
the vertex $1$ from $H$.  Then 
$S(H^{(2)})$ consists of $S(H)\cap G^{(2)}$ together with all
the direct descendants of $1$ in $H$.  This change may be accomplished
by replacing $\mu_b$ by $\rho_b$ whenever $b \in H^{(2)}$ is a 
direct descendant of $1$ in $H$.  Since $1\in S(H)$ for all $H\in \FF''$, 
we obtain the desired equality.
\end{proof}

\begin{example}\label{excont}
We return to Example~\ref{exparred}.  We recall that
the reduced system $\tilde p_1,\dots,\tilde p_6$ is:
\begin{eqnarray*}\label{ex:reduced}
\tilde p_1 = u_1^5  - u_1^3 u_2; \quad
\tilde p_2 = u_2^7 - u_1^8  u_2;\\
\tilde p_3 = u_3^2 - u_1^6 u_3^3;\quad
\tilde p_4 = u_4^3 - u_1^2 u_2 u_5^2 u_6;\\
\tilde p_5 = u_5^2 - u_4^2 u_5 ;\quad
\tilde p_6 = u_6^3 - u_1^3 u_4 u_5 u_6.
\end{eqnarray*}
and, therefore, its associated matrix is
$$B = \left(
\begin{array}{rrrrrr}
2 & -1 & 0 & 0 & 0 & 0\\
-8 & 6 & 0 & 0 & 0 & 0\\
-6 & 0 & -1 & 0 & 0 & 0\\
-2 & -1 & 0 & 3 & -2 & -1\\
0 & 0 & 0 & -2 & 1 & 0\\
-3 & 0 & 0 & -1 & -1 & 2
\end{array}
\right).
$$
Therefore, the system is in normal, triangular form with blocks
relative to the index sets $I_1=\{1,2\}$, $I_2=\{3\}$, and $I_3=\{4,5,6\}$.
The block $B_1$ is global, while $B_2$ and $B_3$ are local.  The graph
$G$ has $3$ vertices $\{1,2,3\}$ and arrows from $1$ to $2$ and 
$1$ to $3$. Hence $S(G) = \{1\}$. The weights are:
$$\delta_1 = 4, \delta_2 = 1, \delta_3 = 5, \rho_1 = 35, 
\rho_2 = 2, \rho_3 = 18\,,$$
and, taking into account the local/global label, we get
$\mu_1 = 31$, $\mu_2 =2$, $\mu_3 = 18$.

%\begin{figure}[ht]
%\begin{center}
%\psfrag{1}{$1$} \psfrag{2}{$2$} \psfrag{3}{$3$}
%\includegraphics[scale=.4]{graph.bmp}
%\caption{The associated graph $G$}
%\end{center}
%\end{figure}

We may now apply (\ref{eq:multat0}) to compute the multiplicity
$\tilde \mu$
of $\langle \tilde p_1,\dots,\tilde p_6\rangle$ at the origin:
$$\tilde \mu = \mu_1 \cdot \rho_2 \cdot \rho_3 = 1116\,.$$

In order to compute $\tilde d$ we use the inductive procedure
of Proposition~\ref{prop:inductive}.  Since the subgraph with
vertices $\{2,3\}$ is disconnected we have:
$$F_2 = (\delta_2 + \mu_2) \cdot (\delta_3 + \mu_3).$$
Hence,  $ \ F_1 = \delta_1\cdot (\delta_2 + \mu_2) \cdot (\delta_3 + \mu_3)
+ \mu_1\cdot \rho_2\cdot \rho_3$.
This gives $\tilde d = 1392$.  We note that this is far from the
B\'ezout bound of $43740$.

Using Lemma~\ref{multred} and Theorem~\ref{parred} we see
that the total number of solutions for the original system
$p_1,\dots,p_8$ are given by $d = 2 \tilde d$ and $\mu = 2 \tilde \mu$. 
This values may be easily verified using a computer algebra
system such as Singular \cite{singularsoft}.

Finally, we note that $G$ has five full subgraphs with vertex sets:
$\{1,2,3\}$, $\{2,3\}$, $\{2\}$, $ \{3\}$,   and $ \emptyset$.  This means that there are
five index sets $\tilde L \subset [6]$, such that $\mu_L \not= 0$.
They are $\tilde L_1 = [6]$, $\tilde L_2 =  \{3,4,5,6\}$, $\tilde L_3 =  \{3\}$, 
$\tilde L_4 = \{4,5,6\}$ and $\tilde L_5 =\emptyset$.  
The corresponding multiplicities are according to (\ref{muL}):
$$\mu_{\tilde L_1} = \tilde \mu = 1116, \ 
\mu_{\tilde L_2} = 144,\ 
\mu_{\tilde L_3} = 40,\ 
\mu_{\tilde L_4} = 72,\ 
\mu_{\tilde L_5} = \tilde\delta = 20.$$
Moreover, the total number of solutions counted  without multiplicity is given
by:
$$\tilde D = \delta_1 + \delta_1 \cdot \delta_3 + \delta_1 \cdot \delta_2 
+ \delta_1 \cdot \delta_2   \cdot \delta_3 = 48$$
This
information may be lifted to the original system using the
bijection $L \to \tilde L$ discussed before Theorem~\ref{parred}. 
 We get that 
$\mu_L=0$ except for the following subsets 
%$L\subset[8]$:
$$L_1 = [8],\ L_2=\{4,5,6,7,8\},\ L_3=\{4\}, \ L_4=\{5,6,7,8\},\ 
L_5=\emptyset.$$
Once again, $\mu_{L_i} = 2\,\mu_{\tilde L_i}$.
\end{example}

\section{Counting complexity}\label{sec:complexity}

In this section we will study the counting complexity,
in the sense of \cite{valiant}, of computing the numerical
invariants $d$, $D$, $\delta$, $\mu$, and $\mu_L$ associated
with a gci $\pn$.  

We have already proved
that we may decide in polynomial time if $\pn$ is a gci
and that the property of being a complete intersection
is independent of the coefficients if $\det B\not=0$.
Moreover, if $\pn$ is a gci we may also transform it into normal and 
triangular form in quadratic time.  Also, since a
system with generic exponents is irreducible and satisfies $\det B\not= 0$, 
we may compute its invariants
in time polynomial in $n$ for any choice of coefficients by Theorem~\ref{vinberg}.  
In the general case,
 we may  compute $\delta$, $\mu$, and $\mu_L$, for a particular choice of $L$, 
directly from the invariants $\delta_a$, $\rho_a$, and $\mu_a$
associated with the diagonal blocks of the system.  Thus, 
$\delta$, $\mu$, and $\mu_L$ may be computed in polynomial
time as well.  

However, we will show below in Theorem~\ref{th:sharpP} that
the computation of $d$ or $D$ is a $\#P$-complete problem, and therefore 
it is at least as hard as an 
NP-complete problem \cite{valiant}.  
In order to do this we begin by 
reversing
the relationship between binomial systems and weighted acyclic directed graphs.  
We recall that to a binomial system $\pn$ in normal and
triangular form we associate an acyclic directed graph $G$ whose
vertices $\{1,\dots,s\}$ correspond to the diagonal blocks of the
associated matrix $B$ and that each vertex has weights 
$\delta_a$, $\rho_a$, $a\in [s]$, plus a label ``local" or ``global".
In the first case we set $\mu_a = \rho_a$, while in the global
case we set $\mu_a = \rho_a - \delta_a$.  In any case $d_a = \delta_a + \mu_a$.  
The proof of the following proposition is straightforward.

\begin{proposition}\label{prop:realizable}
Let  $G= (V,E)$, $V=[s]$, be an acyclic directed graph, with weights
 $\delta_a, \rho_a\in \Z_{>0}$ and labels local/global attached to each vertex. 
Let $\mu_a$ and $d_a$ be defined as above.
 Then, the system of binomials defined by
$$p_a (x_1, \dots, x_s) =  x_a^{d_a} - c_a \left(\prod_{(b,a) \in E} x_b \right) x_a^{\mu_a},$$
for all global vertices $a$, and
$$p_a (x_1, \dots, x_s) =  x_a^{\mu_a} - c_a \left(\prod_{(b,a) \in E} x_b \right) x_a^{d_a},$$
for all local vertices $a$, 
has as weighted graph $(G,\delta_a,\rho_a,\mu_a)$.
\end{proposition}

\begin{remark} The total number of solutions $d$ and $D$ of the system
in Proposition~\ref{prop:realizable} are given by 
(\ref{eq:dformula}) and (\ref{eq:Dformula}), for generic parameters $c_a$.
For any order on the set of
vertices of $G$ such that $i < j$ if there is a path from node
$i$ to node $j$ (i.e., for any linear extension of $G$), it is clear that the corresponding matrix $B$ of
the system will be lower triangular, with diagonal entries $\pm (d_a  - \mu_a)$.
Thus, whenever $d_a \not= \mu_a$, we have that $\det(B) \not=0$ and
we may simply choose $c_a =1$ for all $a \in [s]$.
\par
Note also that 
if $a$ is a source of $G$, then we get 
$p_a = x_a^{d_a} - c_a x_a^{\mu_a}$ in the global case, and
$p_a = x_a^{\mu_a} - c_a x_a^{d_a}$ in the local case. This is
compatible with Remark~\ref{local1}.
\end{remark}

In the particular case when all vertices $\{1,\dots,s\}$ of a directed acyclic
graph $G$ are local, and their weights are $\delta_a=1$, $\rho_a=1$,
for all $a\in [s]$, the binomial system defined in 
Proposition~\ref{prop:realizable} takes a very simple form:
\begin{equation}\label{eq:system}
p_a (x_1, \dots, x_s) =  x_a-  \left(\prod_{(b,a) \in E} x_b \right) x_a^{2}, \  a =1, \dots, s.
\end{equation}
We will refer to this system 
as the {\em standard binomial system} associated with $G$.

\begin{theorem}\label{th:sharpP}
Computing $d$ and $D$ for binomial complete intersections
$\pn$ in normal, triangular form are $\#P$-complete problems.
\end{theorem}

\begin{proof} By Theorem~\ref{th:dformula}, the problems of computing $d$ and
$D$ are in the complexity class $\# P$. We will show that computing these invariants gives,
for special binomial systems, the number of independent subsets
of a bipartite graph $G$.  Since, by \cite{provanball},  this is known to be a $\#P$-complete
problem  the result will follow.

Let $G$ be a bipartite graph with vertices
$\{1,\dots,s\}$. Let $p_1,\dots,p_s$ be the 
standard binomial system of $G$ as in (\ref{eq:system}).
Then, for each full
subgraph $H\subset G$ we have, by (\ref{muL}), that
$\mu_{L(H)} = 1$.  Hence, according to  (\ref{eq:dformula}) and (\ref{eq:Dformula}), 
both $d$ and $D$ are equal to the number of
full subgraphs of $G$.  But, as has been noted earlier, a full subgraph
is completely determined by its sources and, for a bipartite graph $G$,
a subset of vertices is the set of sources of a full subgraph $H$ if and
only if it is an independent subset of $G$.  Thus, $d$ and $D$ agree
with the number of independent subsets of $G$.
\end{proof}

Recall that a directed acyclic graph $G=(V,E)$ 
is called {\em transitive} if  there is an edge
$(a,b) \in E$ each time that there is a directed
path from $a$ to $b$, Transitive directed acyclic graphs are
in correspondence with partial orders $\prec$ on $V$, where $a \prec b$ if and
only if $(a,b) \in E$. Given a partial order $\prec$ on $V$,
a subset $A$ of $V$ is called an antichain if given 
$a_1,a_2\in A$, neither $a_1\prec a_2$, nor $a_2\prec a_1$.
It is shown in \cite{provanball} that counting the number of antichains in
posets is a $\#P$-complete problem and, hence, $\#P$-hard.
Given any directed acyclic graph $G = (V,E)$, it is possible
to compute its transitive closure $G^+ = (V, E^+)$, in time $O(|V|^3)$ by the
well known Floyd--Warshall's algorithm. 
It follows from (\ref{eq:dformula}) and (\ref{eq:Dformula}) that
$d$ and $D$ are the same for the standard binomial systems associated
with $G$ and with $G^+$.

\begin{proposition} 
The 
number of (simple) solutions of the standard system~(\ref{eq:system})
associated with a directed acyclic graph $G$
equals the number of antichains in the associated partial order. 
\end{proposition}

\begin{proof}
As in the proof of Theorem~\ref{th:sharpP}, for the standard binomial
system of $G$ we have $d=D$ and this number agrees with the number
of full subgraphs of $G$.  These subgraphs are determined by their
sources, which correspond exactly to the antichains in the associated
partial order on $V$.
\end{proof}

 Although, as the previous results show, the problem of computing
 the total number of solutions for a general binomial system in
 normal and triangular form is $\#P$-hard, there are classes of
 binomial systems whose invariants may be computed in polynomial
 time.  
For example, if the graph is totally disconnected
 then $d = d_1\cdots d_s = \prod_{i=1}^s (\delta_i + \mu_i)$.  At the
 other extreme if $G$ is a (complete) directed graph with 
vertices $\{1,\dots,s\}$ and  
$(b,a)$ is an edge of $G$ for all $a,b\in [s]$ with $a<b$, then it is easy to 
see that there are only
$s+1$ full subgraphs of $G$ and, consequently, the sums in 
(\ref{eq:Dformula}) and (\ref{eq:dformula}), consist of $s+1$ terms.

Even if the number of full subgraphs is exponential in $s$ and $G$ has
few connected components, a bound on the number of local blocks 
guarantees that $d$ can be computed in polynomial time in $n$. 
For instance, if all blocks are global, then  $B$ is
an $M$-matrix and $\p$ is a Gr\"obner basis for a positive
weight order, and so $d = \rho_1 \cdots\rho_s$.
We end with the following ``positive'' complexity result.

\begin{proposition}\label{polycase}
Let $N \in \Z_{\geq 0}$.
Assume $\p$ is  in normal
and triangular form with $s$ blocks of which at most
$N$ are local. Then, there is a formula to compute 
the total multiplicity $d$ with at most $2^N$ summands, each
involving $s$ products.
Thus,  if the number of local blocks of a binomial system in normal
and triangular form is bounded independently of
$n$, the number of affine solutions of the system can be
computed in time polynomial in $n$.
\end{proposition}

\begin{proof}
Recall the notation 
in Proposition~\ref{prop:inductive}.  
We may write the polynomial formula $F_r ((\delta_a, \rho_a, \mu_a), a \in [s])$ 
for the computation of the
total number of solutions of the system $q^{(r)}$ purely
in terms of $\delta_a$ and $\rho_a$ by keeping track of
the local/global character of each vertex and replacing
$\mu_a$ by $\rho_a$ if $a$ is local and by $\rho_a - \delta_a$
in the case of a global vertex. We call $\tilde{F}_r ((\delta_a, \rho_a), a \in [s])$
the polynomial obtained after these substitutions.
Then, for a global vertex $r$, the recursion~(\ref{eq:Fk}) becomes
\begin{equation} \label{eq:Fkg}
\tilde{F}_r = \delta_r \cdot \tilde{F}_{r+1} + (\rho_r - \delta_r) \cdot
\tilde{F}_{r+1}|_{\delta_a=0},
\end{equation}
where $a$ runs over all direct descendants of $r$.  Let us
write $\tilde{F}_{r+1} = F'_{r+1} \, + F''_{r+1}$, where $F'_{r+1}$
consists of all summands containing a factor $\delta_a$ with
$a$ a direct descendant of $1$. Hence, 
$F'_{r+1}$ vanishes when we set such $\delta_a=0$ and  (\ref{eq:Fkg}) becomes:
$$ \tilde{F}_r = \delta_r \cdot (F'_{r+1} + F''_{r+1}) + (\rho_r - \delta_r) \cdot
F''_{r+1} = \delta_r \cdot F'_{r+1} + \rho_r \cdot
F''_{r+1}$$
and, consequently, the total number of summands does not change 
when adding a global vertex.

On the other hand, if $B_r$ is local then (\ref{eq:Fk}) becomes
\begin{equation*} \label{eq:Fkl}
\tilde{F}_r = \delta_r \cdot \tilde{F}_{r+1} + \rho_r \cdot
\tilde{F}_{r+1}|_{\delta_a=0}
\end{equation*}
and the number of summands is, at worst, doubled.

It follows that when $N$ is bounded independently of the number $n$
of variables, $d$ can be computed by adding  a constant number of summands. Each of these summands has $s \leq n$ products of factors 
involving the computation of determinants  of the square diagonal blocks 
of  the associated matrix $B$ or products of the exponents $r_j$.
\end{proof}

\medskip

\section{Applications}\label{sec:applications}

In this section we will briefly discuss some
of the problems that led us to the study of  systems
of $n$ binomials in $n$ variables.

An important subfamily of binomial ideals is given by
the toric ideals associated to  configurations $A = 
\{a_1,\dots,a_m\}\subset \Z^k$
of integral points spanning $\Z^k$:
$$I_A \, =  \langle x^{u} - x^{v} \ ;\ A\cdot (u-v)=0\rangle\,,$$
where $u,v\in \N^m$.
 In particular, beginning with the work
of Herzog \cite{herzog} and Delorme \cite{delorme}
  the question of classifying complete intersection
toric ideals (and the corresponding semigroup algebras)
has been extensively studied by many authors
\cite{bmt1,  ray, fms1, fms2, fs,  sss}.  A key
step in many of these works is the study of the ideal generated
by binomials $x^{u_i} - x^{v_i}$ associated with a $\Z$-basis
of the kernel of $A$.  More generally, 
given
$\Q$-linearly independent elements 
$\nu_1,\dots,\nu_r \in \Z^m$, consider the associated lattice basis ideal
$J \subset k[x_1,\dots,x_m]$, generated by the binomials
$$b_j = x^{u_j} - x^{v_j}\,;\quad j=1,\dots,r,$$
where $\nu_j = {u_j} - {v_j}$, and ${u_j} ,{v_j}\in \N^m$ have
disjoint support.  Let $\LL \subset \Z^m$ denote the lattice spanned
by $\nu_1,\dots,\nu_r$ and let
$I_\LL \ :=\ \langle x^{u} - x^{v} : u - v \in \LL\rangle$
be the corresponding lattice ideal.  We assume that these ideals
are homogeneous, i.e.,  $w_1+\cdots+w_m=0$, for every $w\in \LL$.

The ideal $I_\LL$ is prime if and only if the lattice
$\LL$ is saturated.  If $\LL$ is not saturated, then 
$I_\LL$ has $g$ radical primary components, where $g$ is
the index of $\LL$ in its saturation.  Moreover, all these
components have the same degree, equal to the degree $d_{\LL}$ of
the associated toric variety \cite{eisenbud}.

We can apply 
Theorem~\ref{th:dformula}      to compute the
multiplicity and geometric degree
\cite{bm} of the primary components of
$J$.  
This may be used to describe the holonomic
rank of Horn systems of hypergeometric partial differential equations and to
study sparse discriminants,
generalizing the  codimension-two case \cite{ ds, dms}.

A straightforward extension of the results of \cite{hs} to
non-saturated lattices gives the following description of all
primary components $\qq$ of $J$.  Let $K\subset \{1,\dots,m\}$ and $Z(K)\subset  \{1,\dots,r\}$
as in (\ref{zeroset}).  Assume that 
$n:=|Z(K)|=|K|$ and for all $j\not \in Z(K)$ 
$${\rm supp}(u_j) \cap K = {\rm supp}(v_j) \cap K = \emptyset.$$
Let $\pp'$  be a primary component of the lattice ideal $I_{\LL'}$
associated to the sublattice
of $\Z^{m-n}$ spanned by $\nu_j$, $j\not \in Z(K)$. 
Then, the ideal 
$$\qq \ = \  \pp' \ + \ \langle b_i, i\in Z(K) \rangle$$
is a primary component of $J$ with associated prime
$$\pp \ =\   \pp' \ + \ \langle  x_k, k\in K \rangle .$$
Note that for $K=\emptyset$ we recover the components 
of $I_{\LL}$.  

In order to describe the multiplicity and geometric degree of
a component $\qq$, let us assume that $K = Z(K) = \{1,\dots,n\}$ and
for any $w\in \Z^m$, denote $\pi(w) = (w_1,\dots,w_n)$.  Let
$\alpha_j = \pi(u_j)$, $\beta_j = \pi(v_j)$ and set 
$$p_j(c;x) \ =\ x^{\alpha_j} - c_j x^{\beta_j}\,, c_j\in k^*.$$
Since $J$ is a complete intersection, $p_1,\dots,p_n$ is a gci.
Let $\mu$ denote the multiplicity at the origin.
Fix coefficients $c\in (k^*)^n$ such that 
$\JJ_c$ is a complete intersection.  
Since 
$$\mu = {\rm length}\left(\kn/\JJ_c\right)_0 = 
{\rm length}\left(k[x_1,\dots,x_m]/J\right)_\pp,$$
and the degree of $\pp$ equals that of $\pp'$, we have

\begin{proposition} 
With notation as above,
the multiplicity of $\qq$ equals $\mu$ and the geometric degree of
$\qq$ equals $d_{\LL'}\cdot \mu$.
\end{proposition}

As a second application, consider a
system of constant coefficient partial differential
equations defined by $n$ operators of the form
\begin{equation}\label{eqsystem}
 a_j \partial^{\alpha_j} - b_j \partial^{\beta_j}\,;\quad  j=1,\dots,n,
\end{equation}
where
$a_j,b_j\in k^* ,\
\alpha_j, \beta_j \in \N^n, \  \alpha_j \not= \beta_j$.
Assume moreover that the ideal $J$ in $\kn$ generated by the binomials
$ a_j x^{\alpha_j} - b_j x^{\beta_j}$ is zero-dimensional.  As before,
let $\mu_L $ the number of points in $\VV(J) \cap \bk^n_L$ counted with
multiplicity.  From
\cite[Chapter 10]{stubook}, we have the following characterization.

\begin{proposition}\label{polysol}
Let $L \subseteq \{1, \dots, n\}$.
The dimension of
the space of solutions to~(\ref{eqsystem}) which depend polynomially
on the variables $x_\ell, \, \ell \in L$, and exponentially on the remaining
variables $x_j, \, j \notin L$, equals $\mu_L$.
\end{proposition}

These dimensions can then be computed using the results in 
Section~\ref{sec:irreducible}, particularly formula (\ref{muL}).

\medskip

\noindent{\bf Acknowledgements: } 
We acknowledge the generous help of many colleagues and friends.
We are grateful to Bernd Sturmfels for the first discussions that lead
to this project. We thank Mircea Mustata for his key suggestions for the
proof of Theorem~\ref{parred}. We are indebted to Peter B\"urgisser
for listening to our questions for many hours, and 
for pointing out the connection of our formulas with the problem of counting
independent sets, which is the key to  our main complexity result.
We also thank  Mart\'{\i}n Mereb and Martin Lotz for useful discussions,
and Daniel Szyld for pointing out the references on non negative
matrices

%\bibliographystyle{plain}

%\bibliography{bci1005}

\end{document}